\documentclass[11pt]{article}
\usepackage{amsmath}
\usepackage{amssymb}
\def\CC{{\rm \kern.24em \vrule width.02em height1.4ex
    depth-.05ex \kern-.26em C}}

\def\TagOnRight

\def\AA{{\it I}\hskip-3pt{\tt A}}

\def\QQ{\rlap {\raise 0.4ex \hbox{$\scriptscriptstyle |$}}
  {\hskip -0.1em Q}}

\newcommand{\be}{\begin{equation}}
\newcommand{\ee}{\end{equation}}
\newcommand{\bea}{\begin{eqnarray}}
\newcommand{\eea}{\end{eqnarray}}
\newcommand{\Bea}{\begin{eqnarray*}}
\newcommand{\Eea}{\end{eqnarray*}}

\catcode`\@=11
\def\theequation{\@arabic{\c@equation}}
\catcode`\@=12

\newcommand{\bi}{\begin{itemize}}
\newcommand{\ei}{\end{itemize}}

\newtheorem{Definition}{Definition}[section]
\newtheorem{Theorem}[Definition]{Theorem}
\newtheorem{Lemma}[Definition]{Lemma}
\newtheorem{Proposition}[Definition]{Proposition}
\newtheorem{Corollary}[Definition]{Corollary}
\newtheorem{Remark}[Definition]{Remark}

\renewcommand{\theequation}{ \thesection.\arabic{equation} }

\begin{document}

\title{\Large Translation Invariant Diffusions in the Space of Tempered Distributions}

\author{B. Rajeev\\
Indian Statistical Institute\\
8th Mile, Mysore Road\\
Bangalore - 560 059\\
email: brajeev@isibang.ac.in}

\maketitle

\begin{abstract} In this paper we prove existence and pathwise uniqueness
for a class of stochastic differential equations (with coefficients
$\sigma_{ij},b_i$ and initial condition $y$ in the space of tempered
distributions) that maybe viewed as a generalisation of Ito's
original equations with smooth coefficients . The solutions are
characterized as the translates of a finite dimensional diffusion
whose coefficients $\sigma_{ij}\star \tilde{y},b_i\star \tilde{y}$
are assumed to be locally Lipshitz.Here $\star$ denotes convolution
and $\tilde{y}$ is the distribution which  on functions, is realised
by the formula $\tilde{y}(r) := y(-r)$ . The expected value of the
solution satisfies a non linear evolution equation which is related
to the forward  Kolmogorov equation associated with the above finite
dimensional diffusion.

\end{abstract}

Processes) 

\section{Introduction}
\setcounter{equation}{0}

In this paper we consider a generalisation of Ito's well known
equation viz. \bea\ dX_t = \sigma (X_t) \cdot dB_t +b(X_t)~dt ~~~~~
;~~~~~  X_0 = x.\eea\ where $\sigma = (\sigma_{ij})_{1 \leq i,j,\leq
d} $ and $ b = (b_1,\cdots,b_d),$ and where $\sigma_{ij}$ and $b_i$
are functions on ${\mathbb R}^d$ and $(B_t)$ is a given
d-dimensional Brownian motion (\cite{KI1}). When $\sigma, b$ are
smooth functions we can use Ito's formula and the duality $\langle
f, \delta_{X_t} \rangle =f(X_t)$ for the random distribution
$\delta_{X_t} \in {\cal S}' $ to arrive at the stochastic
differential equation satisfied by the ${\cal S}'$-valued process
$(Y_t) \equiv (\delta_{X_t})$ viz. \bea dY_t = AY_t\cdot dB_t +
LY_t~dt~~~~;~~~~ Y_0 = y \eea Here $A = (A_1,\cdots A_d) $ , with
$L, A_j$ being  non linear operators from ${\cal S}'$ to ${\cal S}'$
given by \Bea A_j \phi &=& -\sum\limits_{i=1}^d ~ \langle \sigma,
\phi
\rangle_{ij} ~\partial_i ~\phi \\
L \phi &=& \frac{1}{2} ~\sum\limits_{i,j=1}^d~ ( \langle \sigma,\phi
\rangle \langle \sigma, \phi \rangle^t )_{ij}~\partial^2_{ij}~\phi
-\sum\limits_{i=1}^d \langle b, \phi \rangle_i \partial_i \phi \\
\Eea for all $\phi \in {\cal S}'$  where $\langle \sigma,\phi
\rangle_{ij},\langle b, \phi \rangle_i $ is the (i,j)th entry  $
\langle \sigma_{ij},\phi \rangle $ (respectively $i$th entry
$\langle b_i, \phi \rangle $ ) of the matrix $\langle \sigma,\phi
\rangle $ (respectively the vector $\langle b,\phi \rangle $) and
the superscript $t$ denotes transpose. The above equation (1.2)
maybe considered an equation for a system of diffusing particles
,the initial configuration being specified by a tempered
distribution $y $ , with Ito's original equation corresponding to
the position of a single diffusing particle $(X_t)$. For a single
particle ,  the above equation in ${\cal S}'$ is equivalent to Ito's
original equation. This can be seen by action on both sides of (1.2)
by a test function $f \in {\cal S}$ and noting that the resulting
equation is the Ito formula for $f(X_t),$ with $X_t$ given by
equation (1.1). We also note that the Ito formula for $f(X_t)$ is
the point of departure for the celebrated martingale formulation of
Stroock and Varadhan (\cite{SV}). For an alternate approach to
constructing (finite dimensional)`symmetric' diffusions, see
$\cite{F}$

In Section 3, we construct solutions to the above equation in ${\cal
S}'$ when the coefficients $\sigma_{ij},b_i$ and the initial
configuration $y$ are allowed to be elements in ${\cal S}'$ subject
to the condition that the convolutions $\sigma_{ij}\star \tilde{y},
b_i\star \tilde{y}$ be locally Lipshitz  functions on ${\mathbb R}^d
$(see Theorem 3.4, Remark 3.7 below). We may take $\sigma_{ij}$ and
$b_i$ to be arbitrary distributions - for example a Dirac
distribution - provided $y$ is sufficiently smooth or conversely $y$
to be `arbitrary' and $\sigma_{ij},b_i$ to be smooth. Our methods
also extend to the case when the coefficients of $L$ and $A_i$ are
non linear functions (see Remark 3.9).  We use the countable
Hilbertian structure of ${\cal S}'$ viz. the fact that ${\cal S}'  =
\bigcup\limits_{p \in {\mathbb R}} {\cal S}_p$, where the ${\cal
S}_p$ are Hilbert spaces, to formulate our equations
($\cite{KI2},\cite{KX}$).The solutions $(Y_t(y))$ of equation (1.2),
starting at $y$ at time zero , are given explicitly as follows : Let
$\tau_x , x \in {\mathcal R}^d$ denote the translation operators
acting on ${\cal S}'$. Then $Y_t(y) = \tau_{z_t(y)}(y)$ where
$(z_t(y))$ is a finite d- dimensional diffusion starting at zero and
whose drift and diffusion coefficients are given by the convolutions
viz. $\sigma_{ij}\star \tilde{y}$ and $b_i\star \tilde{y}$. Note
that since $z_t(y)$ depends on $y$, the translation
$\tau_{z_t(y)}(y)$ is non linear (unlike the Brownian case
\cite{RT1}). Existence and pathwise uniqueness of $(Y_t(y))$are then
a consequence of existence and pathwise uniqueness for $(z_t(y))$
and the Ito formula ( $\cite{BR}$). Our proof of the above
representation uses the uniqueness of solutions of equation (1.2)
with random operator coefficients
$\bar{A}(s,\omega),\bar{L}(s,\omega)$ (see $\cite{P},\cite{KR}$  )
that satisfy the monotonicity inequality (see $\cite{GMR}$ and
Theorem 3.5 below).

In Corollary(3.8), we show that $Y_t(y)$ has the translational
invariance property viz. the solution corresponding to an initial
condition which is a translate of $y \in {\cal S}^{\prime}$ by $x
\in {\mathbb R}^d$ viz. $Y_t(\tau_x(y)) $ is equal to $\tau_{x
+z_t(\tau_x(y))}(y)$ and $X_t(x,y) :=  x +z_t(\tau_x(y)) $ is the
solution of the SDE starting at $x$ at $t=0,$ which is satisfied by
$(z_t(y))$ for $x=0.$ Consequently, the state space ${\cal S}_p$ of
the flow generated by the solutions of (1.2) viz. $(Y_t(y))$ splits
up into `irreducible' components $C(y) := \{\tau_x(y); x \in
{\mathbb R}^d \}$ on each of which $(Y_t(y))$ behaves as the finite
dimensional diffusion $(X_t(x,y)).$ This allows us to describe the
Markov properties of $(Y_t(y))$ in terms of that of that of
$(X_t(x,y)).$ We do this in Section 4.

In Section 5, Theorem 5.4, we study the deterministic non linear
evolution equation that results on taking expectations in equation
(1.2). We represent the solutions  in terms of a `non linear
convolution'. These are related to the fundamental solutions of the
forward equation associated with the diffusion $X_t(x,y)$ (see
Remark 5.5) We describe the latter in Theorem 5.6. This extends some
results of $\cite{RT2}$ to the case when the coefficients are non
smooth.

Physically, the form of the solutions viz. $Y_t = \tau_{z_t}y$
suggests a conservation law. Indeed, when the initial distribution
$y$ is given by an integrable function $y(r)$ and when the solution
$Y_t$ or equivalently $z_t$ is defined for all $t \geq 0$ i.e. there
is no explosion, then we have for all $t \geq 0$, $$ \int_{{\mathbb
R}^d}Y_t(r)~dr = \int_{{\mathbb R}^d}y(r)~dr .$$  The form of the
solution further suggests that the diffusion of a system of
particles according to (1.2), starting in an initial configuration
given say, by a function $y = y(r)$, in a diffusive  medium
represented by the coefficients $\sigma_{ij}$ and $b_i$ is
equivalent to a translational flow of the same initial
configuration, in a new medium -  represented by coefficients
$\sigma_{ij}\star \tilde{y}$ and $b_i\star \tilde{y}$ - which is the
result of an interaction between the initial particles and the
original diffusive medium. We refer to \cite{VSB} for some related
ideas.

\section{Preliminaries}
Let $\left(\Omega , {\mathcal F}, \left\{{\mathcal
F}_t\right\}_{t\geq 0}, P\right)$ be a filtered probability  space
satisfying the {\it usual conditions} viz. 1) $\left(\Omega ,
{\mathcal F},  P\right)$ is a complete probability space. 2)
${\mathcal F}_0$ contains all $A\in {\mathcal F}$, such that $P(A) =
0$, and 3) ${\mathcal F}_t = {\displaystyle \bigcap_{s>t}}{\mathcal
F}_s$. On this probability space is given a standard d-dimensional
${\cal F}_t$- Brownian motion $(B_t) \equiv (B^1_t, \ldots, B_t^d)$.
We will denote the filtration generated by $(B_t)$ as $({\cal
F}_t^B).$ Consider now the equation (1.1). Let
$\bar{\sigma}_{ij},\bar{b}_i$ be locally Lipshitz functions on
${\mathbb R}^d$ for $i,j = 1,\cdots,d$. Let $\bar{\sigma} :=
(\bar{\sigma}_{ij})$ and $(\bar{b}) := (\bar{b}_1,\cdots,\bar{b}_d)$
. We use the notation $\hat {\mathbb R}^d = {\mathbb R}^d \cup
\{\infty \}$ for the one point compactification of ${\mathbb R}^d.$

\begin{Theorem} Let $\bar{\sigma}, \bar{b}, (B_t)$ be as above. Then
$\exists~ \eta: \Omega \rightarrow (0,\infty], \eta$ an $({\cal
F}_t^B)$ stopping time and an $\hat {\mathbb R}^d$-valued,
$({\cal F}^B_t)$ adapted process $(X_t)_{t \geq 0}$ such that \\
\begin{enumerate}
\item For all $\omega \in \Omega$,~ $X_.(\omega) : [0,\eta(\omega)) \rightarrow {\mathbb R}^d,
$~is continuous
 and $X_t(\omega)
= \infty,~~ t \geq \eta(\omega) $\\
\item a.s. (P), $\eta(\omega) < \infty $ implies $\lim\limits_{t \uparrow \eta(\omega)} X_t(\omega) = \infty
.$\\
\item a.s.(P),\bea  X_t = x+\int\limits_0^t \bar{\sigma}  (X_s) \cdot dB_s
+\int\limits_0^t \bar{b}(X_s)~ds \eea for $0 \leq t < \eta
(\omega)$.
\end{enumerate}
 The solution $(X_t, \eta)$ is (pathwise)
unique i.e. if $(X^1_t,\eta^1)$ is another solution then $P\{ \eta =
\eta_1 , X_t = X^1_t, 0 \leq t < \eta\} = 1.$
\end{Theorem}
{\bf Proof :} We refer to \cite{IW}, Chapter IV, for the proofs.
Theorem 2.3 for the proof of existence and Theorem 3.1 for the proof
of uniqueness .$\hfill{\Box}$

Let $\alpha,\beta  \in {\mathbb Z}^d_+ := \{(x_1,\cdots,x_d) : x_i
\geq 0, x_i \rm {integer}\}.$ Let  $x^\alpha := x_1^{\alpha_1}
\ldots x_d^{\alpha_d}$ and $\partial^\beta :=\partial_1^{\beta_1}
\ldots
\partial_d^{\beta_d} .$ For a multi index $\alpha$, we use the notation $|\alpha| :=
\sum\limits_{i=1}^d \alpha_i.$ Let ${\mathcal S}$ denote the space
of rapidly decreasing smooth real functions on ${\mathbb R}^d$ with
the topology given by the family of semi norms
$\wedge_{\alpha,\beta}$,defined for $f \in {\mathcal S}$ and multi
indices $\alpha, \beta$ by $\wedge_{\alpha,\beta} (f) :=
\sup\limits_x ~|x^\alpha
\partial^\beta  f(x) |$. Then $\{{\cal S}, \wedge_{\alpha,\beta}\}: \alpha, \beta \in
{\mathbb Z}^d_+ \}$ is a Fr\'{e}chet space. ${\mathcal S}'$ will
denote its continuous dual. The duality between ${\cal S}$ and
${\cal S}'$ will be denoted by $\langle \psi , \phi \rangle$ for
$\phi \in {\cal S}$ and $\psi \in {\cal S}'$. For $x \in {\mathbb
R}^d $ the translation operators $\tau_x : {\cal S} \rightarrow
{\cal S} $ are defined as $ \tau_x f(y) :=f(y-x)$ for $f \in {\cal
S}$ and then for $\phi \in {\cal S}'$ by duality : $ \langle\tau_x
\phi,f\rangle := \langle \phi,\tau_{-x} f \rangle $.Let $\{h_k; k
\in {\mathbb Z}^d_+ \}$ be the orthonormal basis in the real Hilbert
space  $L^2 ({\mathbb R}^d,dx) \supset {\cal S}$ consisting of the
Hermite functions (see for eg. \cite{T}); here $dx$ denotes Lebesgue
measure. Let $\langle \cdot,\cdot \rangle_0$ be  the inner product
in $L^2 ({\mathbb R}^d,dx)$. For $f \in {\cal S}$ and $p \in
{\mathbb R}$ define the inner product $\langle f, g\rangle_p $ on
${\cal S}$ as follows : \Bea \langle f, g\rangle_p
:=\sum\limits_{k=(k_1,\cdots,k_d)\in {\mathbb Z}^d_+}(2 |k|+d)^{2p}
~\langle f, h_k \rangle_0 ~\langle g, h_k \rangle_0 \Eea The
corresponding norm will be denoted by  $ \|\cdot \|_p$. We define
the Hilbert space ${\cal S}_p$ as the completion of $ {\cal S}$ with
respect to the norm  $ \|\cdot \|_p$. The following basic relations
hold between the ${\cal S}_p$ spaces (see for eg. \cite{KI2},
\cite{KX}): For $0 < q < p, ~{\cal S} \subset {\cal S}_p \subset
{\cal S}_q \subset L^2 = {\cal S}_0 \subset {\cal S}_{-q} \subset
{\cal S}_{-p} \subset {\cal S}'$. Further,
 ${\cal S}' =
\bigcup\limits_{p \in {\mathbb R}} {\cal S}_p$ and
$\bigcap\limits_{p \in {\mathbb R}} {\cal S}_p ={\cal S}$. If
$\{h_k^p: k \in {\mathbb Z}^d_+ \}$ denotes the orthonormal basis in
${\cal S}_p$ consisting of the (suitably normalised) Hermite
functions, then the dual space ${\cal S}'_p$ maybe identified with
${\cal S}_{-p}$ , via the  basis $\{h_k^{-p}: k \in {\mathbb Z}^d_+
\}$ of ${\cal S}_{-p}$. For $\phi \in {\cal S}$ and $\psi \in {\cal
S}'$ the bilinear form $(\phi,\psi) \rightarrow \langle \psi,
\phi\rangle$ also gives the duality between ${\cal S}_p (\supset
{\cal S})$ and ${\cal S}_{-p} (\subset {\cal S}')$. It is also well
known that $\partial_i : {\cal S}_p \rightarrow {\cal
S}_{p-\frac{1}{2}}$ are bounded linear operators for every $p \in
{\mathbb R}$ and $i = 1,\cdots,d$. Suppose for $i,j = 1,\cdots,d$,
$\sigma_{ij},b_i \in {\cal S}_{-p}$. Let $\sigma :=
(\sigma_{ij})_{1,\leq i,j,\leq d}$ and $b := (b_1 \cdots b_d).$We
can then define for $j= 1,\cdots,d$ the non linear operators $A_j,
L: S_p \rightarrow S_{p-1}$ as follows: \Bea A_j \phi &:=&
-\sum\limits^d_{i=1}~ \langle \sigma, \phi\rangle_{ij}
~\partial_i \phi \\
L\phi & :=& \frac{1}{2}~\sum\limits^d_{i,j=1} (\langle \sigma, \phi
\rangle \langle \sigma, \phi \rangle^t)_{ij}~ \partial^2_{ij} \phi~
-\sum\limits^d_{i=1}~\langle b, \phi \rangle_i~ \partial_i \phi.\Eea
Here  $\langle \sigma, \phi \rangle $ is the matrix whose $(i,j)$th
entry $\langle \sigma, \phi\rangle_{ij} $ is $\langle \sigma_{ij},
\phi \rangle $ , $\langle \sigma, \phi \rangle^t $  its transpose
and $\langle b,\phi \rangle $ is the vector whose $i$th entry
$\langle b,\phi \rangle_i$ is $\langle b_i,\phi \rangle $. The d-
tuple of operators $(A_1 , \ldots , A_d)$ will denoted by $A$.

\section{Stochastic Differential Equations in ${\cal S'}$}
 We now consider a stochastic partial differential equation in
${\cal S}'$ driven by the Brownian motion $(B_t)$ and differential
operators $A_i,i = 1,\cdots,d$ and $L$ defined above with given
coefficients $\sigma_{ij},b_i, i,j = 1,\cdots d$ in the space ${\cal
S}_{-p}$ for some fixed $p \in {\mathbb R}$ and initial condition $y
\in S_p$ viz. \bea dY_t = A(Y_t) \cdot dB_t +L(Y_t)~dt ~~~~;~~~~ Y_0
=y. \eea  Note that if $(Y_t)$ is an ${\cal S}_p$ valued ,locally
bounded , $({\cal F}_t)$ adapted process then
$A_i(Y_s),i=1,\cdots,d$ and $L(Y_s)$ are ${\cal S}_{p-1}$ valued,
adapted , locally bounded processes and hence the stochastic
integrals $\int_0^tA_i(Y_s)~dB^i_s$ and $\int_0^tL(Y_s)~ds$ are well
defined ${\cal S}_{p-1}$ valued, continuous ${\cal F}_t$ adapted
processes and in addition, the former processes are ${\cal F}_t$
local martingales . We then have the following definition of a
`local' strong solution of equation(3.4) .

\begin{Definition}\rm{ Let $p \in {\mathbb R}.$ Let $y \in {\cal S}_p$, $\sigma_{ij}, b_i
\in {\cal S}_{-p}$ and $\{B_t,{\cal F}_t\}$ the given standard
$({\cal F}_t)$ Brownian motion. Let $\delta$ be an arbitrary
state,viewed as an isolated point of  $\hat{{\cal S}_p} := {\cal
S}_p \cup \{\delta\}.$ By an $\hat{{\cal S}_p}$ valued, strong,
(local) solution of equation (3.4), we mean a pair $(Y_t(y), \eta)$
where $\eta : \Omega \rightarrow (0,\infty]$ is an ${\cal F}_t^B$
stop time and $(Y_t(y))$ an
 $\hat{{\cal S}_p} $ valued $({\cal F}_t^B)$
 adapted process such that \begin{enumerate} \item  For all $\omega \in \Omega
 ,$ $Y_.(y,\omega) : [0,\eta(\omega)) \rightarrow {\cal S}_p$ is a continuous map
  and
 $Y_t(y,\omega) = \delta , t \geq \eta(\omega)$
 \item a.s. (P)\rm{ the following equation holds in ${\cal S}_{p-1}$ for $0
\leq t < \eta (\omega)$}, \bea Y_t(y) = y +
\sum\limits_{j=1}^d\int\limits_0^t A_j(Y_s(y)) ~dB_s^j +
\int\limits_0^t L(Y_s(y))~ds. \eea
\end{enumerate}}
\end{Definition}

\begin{Remark} \rm {By usual stopping arguments, the $S_{p-1}$ valued
stochastic integrals on the right hand side in equation

(3.5) can be shown to be finite , almost surely , for $0 \leq t <
\eta.$ Also , all functions $f:{\cal S}_p \rightarrow {\mathbb R}$
are extended to functions on $\hat{{\cal S}_p} $ by setting
$f(\delta) := 0.$}
\end{Remark}

\begin{Definition}
{\rm We say that pathwise uniqueness holds for equation (3.4)  iff
given $(B_t)$, a standard Brownian motion on a probability space
$(\Omega,{\cal F},P)$ , and given $y \in {\cal S}_p$, $\sigma_{ij},
b_i \in {\cal S}_{-p}$, and any two $S_p$ valued strong solutions
$(Y^i_t(y),\eta^i), i = 1,2$ of equation (3.4), we have  $P\{
Y^1_s(y) = Y^2_s(y), 0 \leq  s < \eta^1\wedge\eta^2 \} = 1.$}
\end{Definition}

\begin{Theorem} {Let $y \in {\cal S}_p$, $\sigma_{ij}, b_i
\in {\cal S}_{-p}$ and $\{B_t,{\cal F}_t\}$ the given standard
$({\cal F}_t)$ Brownian motion. Suppose the ${\mathbb R}^{d^2}$
valued function $\bar{\sigma}(x) := (\langle \sigma_{ij},\tau_xy
\rangle)$ and the ${\mathbb R}^d$ valued function $\bar{b}(x)=
(\langle b_i ,\tau_xy \rangle )$ are locally Lipshitz. Then equation
(3.4) has a unique strong solution.}

\end{Theorem}

To prove  the uniqueness assertion in Theorem 3.4 , we need Lemma
(3.6) below . It is of independent interest since it characterises
the solutions of equation (3.4). The proof of the lemma in turn
depends on the so called `Monotonicity inequality' which we now
state. Given real numbers $\sigma_{ij},b_i, i,j = 1,\cdots d,$ let
$\sigma = (\sigma_{ij})$ and $\sigma^t $ the transpose of $\sigma.$
We define the constant coefficient differential operators $A_{oj},
j= 1 \cdots d$ and $L_o$ as follows : \Bea A_{oj} \phi
=-\sum\limits^d_{i=1}~ \sigma_{ij} ~\partial_i \phi, ~~j=1,\ldots,d
\Eea \Bea  L_o \phi =\frac{1}{2} ~\sum\limits^d_{i,j=1}~ (\sigma
\sigma^t)_{ij}~
\partial^2_{ij} \phi -\sum\limits^d_{i=1}  b_i ~\partial_i \phi. \Eea

\begin{Theorem}
{ Let $\alpha > 0.$Then there exists a constant $C = C(\alpha,p,d)
>0$ depending only on $\alpha$,$p$ and $d$ such that   for all $\phi \in {\cal S}_p$,
and all $\sigma_{i,j},b_i, i,j = 1,\cdots, d$ with  $\alpha \geq
\max\limits_{i,j}\{ |\sigma_{i,j}|,|b_i|\}, $ we have \Bea 2 \langle
\phi, {L_o} \phi \rangle_{p-1} +\sum\limits^d_{i=1}\|{A}_{oi} \phi
\|^2_{p-1} \leq C ~\|\phi\|^2_{p-1} \Eea}
\end{Theorem}

{\bf Proof :} See \cite{GMR}.

\begin{Lemma}  {Suppose (3.4) has an ${\cal S}_p$-valued
strong solution $(Y_t(y),\eta)$. Define the continuous
semi-martingale in ${\mathbb R}^d$ as follows: \Bea   z_t(y) :=
\int\limits_0^{t\wedge \eta} \langle \sigma, Y_s \rangle \cdot dB_s
+\int\limits_0^{t\wedge \eta} \langle b, Y_s \rangle ~ds. \Eea Then
a.s P,~ $Y_t(y) = \tau_{z_t(y)}(y) $, for  $0 \leq t < \eta$ .}
\end{Lemma}

{\bf Proof:} We use the notation $Y_t , z_t$ for $Y_t(y),z_t(y)$
respectively. We apply the Ito formula in Theorem 2.3 of $\cite{BR}$
to the process $(\tau_{z_t}y)$ to get the following equation in
${\cal S}_{p-1}$ : \Bea \tau_{z_t}y &=& \tau_{z_0}y
-\int\limits_0^t \partial_i (\tau_{z_s}y) ~dz^i_s \\
&& + \frac{1}{2}~ \sum\limits_{i,j =1}^d ~\int\limits_0^t
\partial^2_{ij}~(\tau_{z_s}y) ~d\langle z^i, z^j \rangle_s  \\
&=& y + \int\limits_0^t \bar A (s,\omega) (\tau_{z_s}y) \cdot
dB_s \\
&& + \int\limits_0^t \bar L (s,\omega) (\tau_{z_s}y) ~ds \Eea where
$(\bar{A}_j(s)), j = 1 \cdots d$ and $(\bar{L}(s))$ are operator
valued ${\cal F}_t$ adapted processes such that for $j=1,\cdots,d$
and each $(s,\omega)$, $\bar A_j(s,\omega)$and $ \bar L(s,\omega)$
are linear operators from $ {\cal S}_{p}$ into ${\cal S}_{p-1}$ ,
defined as follows : \Bea \bar A_j (s,\omega) \phi
=-\sum\limits^d_{i=1}~\langle \sigma, Y_s(\omega) \rangle_{ij}
~\partial_i \phi, ~~j=1,\ldots,d \Eea \Bea \bar L(s,\omega) \phi
=\frac{1}{2} ~\sum\limits^d_{i,j=1}~(\langle
\sigma,Y_s(\omega)\rangle \langle \sigma, Y_s(\omega) \rangle^t
)_{ij}~
\partial^2_{ij} \phi -\sum\limits^d_{i=1} \langle b, Y_s (\omega)
\rangle_i ~\partial_i \phi. \Eea Note that almost surely for $0 \leq
t < \eta$, $(Y_t)$ satisfies the same equation as $(\tau_{z_t}y)$ in
${\cal S}_{p-1}$ i.e.  \Bea Y_t = y+\int\limits_0^t \bar A (s)Y_s
\cdot dB_s +\int\limits_0^t \bar L(s)Y_s ~ds. \Eea Let $(\sigma_n)$
be a sequence of $({\cal F}_t)$ stopping times, such that almost
surely, $0 \leq \sigma_n < \eta$ , $\sigma_n \uparrow \eta$ and such
that for each $n \geq 1,$ the processes
$(\langle\sigma_{ij},Y_{\sigma_n\wedge s}\rangle)$and $(\langle b_i
,Y_{\sigma_n\wedge s}\rangle)$ are uniformly bounded for all $i,j$.
Define for $n \geq 1,$ $X^n_t := Y_{\sigma_n\wedge t} -
\tau_{z_{\sigma_n\wedge t}}y$. Let for $n \geq 1,$$C_n$ be the
constant obtained from Theorem 3.5 above applied to the operators
$A_{oi} = A_i(s,\omega), i = 1,\cdots,d $ and $ L_o = L(s,\omega)$
for fixed $(s,\omega)$ with $s \leq \sigma_n(\omega), n \geq 1$ and
$$\alpha_n
> \sup\limits_{i,j, s ,\omega} \{\langle \sigma_{ij}, Y_s(\omega) \rangle,
 \langle b_i, Y_s (\omega)
\rangle\}$$ where the supremum is taken over $1 \leq i,j,\leq d$, $
s \leq \sigma_n(\omega)$ and $ \omega \in \Omega$.We then have,
using integration by parts,

 \Bea \|X^n_{t \wedge \sigma_n}\|^2_{p-1} &=& \int\limits_0^{t \wedge
\sigma_n} \left( \sum\limits_{i=1}^d \|\bar{A}_i (s) X^n_s\|_{p-1}
+2 \langle X^n_s , \bar{L}(s)~X^n_s \rangle_{p-1} \right) ~ds
+M^n_{t \wedge \sigma_n}\\
&\leq& C_n ~\int\limits_0^{t \wedge \sigma_n} \|X^n_s \|^2_{p-1} ~ds
+M^n_{t \wedge \sigma_n} \\
&\leq& C_n ~\int\limits_0^t \|X^n_{s \wedge \sigma_n}\|^2_{p-1}~ds
+M^n_{t \wedge \sigma_n} \\ \Eea where $(M^n_t)$ is a continuous
local martingale.  The Gronwall inequality now implies that for
every $t \geq 0 $ and $n \geq 1$, $\Rightarrow E\|X^n_{t \wedge
\sigma_n}\|^2_{p-1}=0 $. It follows that $ Y_t = \tau_{z_t}(y)$
almost surely for $0 \leq t < \eta $.$\hfill{\Box}$

\begin{Remark}\rm { For a function $f$ we use the notation $\tilde{f}$ to
denote the function $\tilde{f}(r) := f(-r)$. For a distribution $y$
we define the distribution $\tilde{y}$ via duality i.e. $\langle
\tilde{y}, f \rangle := \langle y , \tilde{f}\rangle ~~ \forall f
\in {\cal S}.$ We note that $\bar{\sigma_{ij}}(x)  := \langle
\sigma_{ij},\tau_xy \rangle = \sigma_{ij}\star \tilde{y}(x)$ where
$\star$ denotes convolution . It is well known that when
$\sigma_{ij} \in {\cal S}'$ and $y \in {\cal S}$ the convolution
$\sigma_{ij}\star y$ is ${\cal C}^{\infty}$(see Theorem 30.2,
\cite{FT}).}
\end{Remark}
We have the following Corollary to Lemma 3.6.

\begin{Corollary} Let $y \in {\cal S}_p, \sigma_{ij},b_i \in {\cal
S}_{-p}, i,j = 1,\cdots ,d $ and $ x \in {\mathbb R}^d.$ Suppose
that the coefficients in equation (2.3) are given as $\bar
\sigma_{ij} = \sigma_{ij}\star \tilde{y}$ and $\bar b_{i} =
b_{i}\star \tilde{y} $. Let $(Y_t(\tau_x(y)),\eta)$ be  a solution
of equation (3.4) with initial condition $\tau_x(y).$ Let
$(z_t(\tau_{x}(y)))$ be given by Lemma 3.6  with $ Y_t =
\tau_{z_t(\tau_x(y))}(\tau_x(y)) = \tau_{x+z_t(\tau_x(y))}(y), t <
\eta $. Then the process $ (x + z_t(\tau_x(y)))$ solves equation
(2.3) upto the random time $\eta$, with initial condition $x $.

\end{Corollary}

{\bf Proof } Let $X_t := x + z_t(\tau_x(y)), t < \eta .$ We have
a.s. for $t < \eta$, \Bea X_t = x + z_t(\tau_x(y)) &=& x +
\int\limits_0^t \langle \sigma, \tau_{ z_s(\tau_x(y))}(\tau_x(y)
\rangle \cdot dB_s +\int\limits_0^t \langle b,
\tau_{z_s(\tau_x(y))}(\tau_x(y)) \rangle~ds \\ &= & x +
\int\limits_0^t \langle \sigma, \tau_{x+ z_s(\tau_x(y))}(y) \rangle
\cdot dB_s +\int\limits_0^t \langle b, \tau_{x+z_s(\tau_x(y))}(y)
\rangle~ds
\\ &=& x + \int\limits_0^t \bar \sigma (X_s) \cdot dB_s +\int\limits_0^t
\bar b (X_s) ~ds \Eea  Since $\bar \sigma (x) =\langle \sigma,
\tau_x(y) \rangle$ and $ \bar b =\langle b, \tau_x(y) \rangle$, the
result follows from Theorem 2.1.$\hfill{\Box}$

{\bf Proof of Theorem 3.4} Assume that $\bar \sigma_{ij}$ and $\bar
b_i$ are locally Lipschitz functions on ${\mathbb R}^d.$.
 Let $(z_t(y),\eta )$  be a solution of equation (2.3). $z_t := z_t(y).$
 Define the $S_p$ valued process
 $(Y_t)$ as follows : \Bea Y_t := \tau_{z_t}y, && t < \eta \\
\delta ~~~~ && t \geq \eta. \Eea Note that $\bar {\sigma}_{ij} (z_t)
=\langle {\sigma}_{ij}, \tau_{z_t}y \rangle =\langle \sigma_{ij},
Y_t\rangle$ and similarly  $\bar {b}_i(z_t) =\langle b_i,
Y_t\rangle, t < \eta.$ Applying the Ito formula in Theorem 2.3 of
$\cite{BR}$ to the process $(\tau_{z_t}y)$ as in the proof of Lemma
3.6, we get , almost surely for $t < \eta$,

 \Bea Y_t =\tau_{z_t}y &=&\tau_{z_0}y - \sum\limits_{i=1}^d \int\limits_0^t
\partial_i \tau_{z_s}y ~dz_s^i +\frac{1}{2}
~\sum\limits_{i,j=1}^d~\int\limits_0^t
\partial^2_{ij} \tau_{z_s}y ~d\langle z^i, z^j \rangle_s.\\
&=& y + \sum\limits_{j=1}^d\int\limits_0^t A_j(Y_s) ~dB_s^j +
\int\limits_0^t L(Y_s)~ds.
 \Eea

Hence $(Y_t ,\eta)$ is a solution. Suppose $(Y^1_t, \eta^1)$ and
$(Y^2_t, \eta^2)$ are two solutions. Define  the processes $(z^i_t)$
as follows : For $0 \leq t < \eta^i,$\Bea z^i_t
:=\int\limits_0^{t\wedge \eta^i} \langle \sigma, Y^i_s \rangle \cdot
dB_s^i &+& \int\limits_0^{t\wedge \eta^i} \langle b,Y^i_s\rangle ~ds
. \Eea We take $z^i_t := \infty $ for $t > \eta^i.$Then by Lemma 3.6
,$Y^i (t) =\tau_{z^i_t}(\phi), t < \eta^i.$ But then by Cor (3.8),$
(z^i_t, \eta^i), i=1,2$ solves equation (2.3) with initial condition
$x =0,$ upto $ \eta^i, i = 1,2.$ The proof of uniqueness in Theorem
(3.1) of \cite{IW} can be applied here (with appropriate
modification to take care of the limit at $\eta^1\wedge\eta^2$) to
conclude that, almost surely , $z^1_t = z^2_t $ for $0\leq t <
\eta^1\wedge\eta^2 $. It follows that almost surely ,$ Y^1_t = Y^2_t
$ for $~0\leq  t < \eta^1\wedge\eta^2 $.$\hfill{\Box}$
\begin{Remark}\rm{Instead of the linear functionals $\langle \sigma_{ij} , .
 \rangle, \langle b_i,.\rangle :{\cal S}_p
\rightarrow {\mathbb R}$ as coefficients of $L$ and $A_i$, we could
take more general (non linear) functions $\sigma_{ij},~b_i :{\cal
S}_p \rightarrow {\mathbb R}$ as the coefficients in $L$ and the
$A_i$'s . Theorem 3.4 remains true provided $\sigma_{ij}(\tau_x y),
b_i(\tau_xy)$ are locally Lipschitz functions on ${\mathbb R}^d$ for
$y \in {\cal S}_p$.}
\end{Remark}

Consider the  ordinary differential equation $$ \dot{y} = b(y) ~;
~~~ y(0) = y_o$$ where $b = (b_1,\cdots,b_d)$ is a smooth vector
field on ${\mathbb R}^d.$ The following corollary generalises this
equation when the $b_i$ are tempered distributions on ${\mathbb
R}^d.$
\begin{Corollary} {Suppose $\sigma \equiv 0, b_i \in S_{-p},y \in S_p$.
 Then the unique solution of the  first order evolution equation viz.
\Bea \partial_t Y_t &=& -\sum\limits_{i=1}^d \langle b_i, Y_t
\rangle~
\partial_i Y_t,; ~~Y_0 = y. \\
&=& -\langle b,Y_t \rangle \cdot \nabla Y_t,; ~~ Y_0 = y.  \Eea is
given by $\{Y_t, 0 \leq t < \eta \} $ where $ Y_t := \tau_{z_t}y $
for $0 \leq t < \eta$ and  $\{z_t , 0 \leq t < \eta\}$ is the unique
solution of the ordinary differential equation
$$
 \dot z_t  =\bar b (z_t), ~~ z_0 = 0. $$}
\end{Corollary}

We can characterise the unique global solutions of equation (3.4) in
terms of the unique solutions of equation (2.3) described in Theorem
2.1.
\begin{Theorem} Let $\sigma_{ij},b_i, \bar{\sigma}_{ij},\bar{b}_i,y$
 be as in Theorem 3.4. Then
there exists a solution $(Y_t(y),\eta)$ of equation (3.4) with the
property that $Y_t(y) = \tau_{z_t(y)}(y), t < \eta,$ where
$(z_t(y),\eta)$ is the unique solution of equation (2.3) given by
Theorem (2.1). In particular, almost surely, on the set $\{\eta <
\infty\},~\lim\limits_{t \uparrow \eta }z_t(y) = \infty.$ The
solution is pathwise unique, i.e. if $(Y^1_t,\eta^1)$ is another
such solution, then $P\{ \eta = \eta^1, Y_t(y) = Y^1_t, 0 \leq t <
\eta \} = 1.$
\end{Theorem}
{\bf Proof :} Let $(z_t(y),\eta)$ be the unique solution of
equation(2.3) given by Theorem(2.1) and define $Y_t(y) :=
\tau_{z_t}(y), t < \eta$ and $Y_t(y) := \delta , t \geq \eta.$ Then
as in the proof of Theorem(3.4), $(Y_t(y),\eta)$ is a solution of
equation(3.4). The uniqueness follows from the uniqueness of
$(z_t(y),\eta).$ $\hfill{\Box}$

 The following proposition charactarises  the global solutions of equation(3.4)
 as processes in ${\cal S}^{\prime}$  . It shows that they are
 better behaved than their finite dimensional counterparts. It
 concerns the purely analytical behaviour of
 $\tau_xy$ as $|x| \rightarrow \infty.$

\begin{Proposition}
Let $\sigma_{i,j},b_i,y$ be as in Theorem 3.4, with $y \neq 0$. Let
$Y_t(y) = \tau_{z_t(y)}y$ be the unique solution upto time $\eta$ of
equation(3.4), where $(z_t(y))$ solves equation(2.3) with $x = 0$
and $\bar \sigma_{ij},\bar b_{i} $ as in Corollary (3.8). Fix $
\omega \in \Omega.$ Then, $z_t(y,\omega) \rightarrow \infty $ as $t
\rightarrow \eta(\omega)$ whenever $Y_t(y,\omega) \rightarrow 0$
weakly in ${\cal S}^{\prime}$ as $t \rightarrow \eta(\omega).$
Conversely, suppose one of the following two conditions is satisfied
viz. \begin{enumerate} \item $y$ is a square integrable function
i.e. $y \in L^2({\mathbb R}^d)= {\cal S}_0.$
\item $y$ has compact support i.e. $y \in {\cal E}^\prime.$ \end{enumerate}
Then, $z_t(y,\omega) \rightarrow \infty $ as $t \rightarrow
\eta(\omega)$ implies $Y_t(y,\omega) \rightarrow 0$ weakly in ${\cal
S}^{\prime}.$
\end{Proposition}
{\bf Proof :} Let $Y_t(\omega):= Y_t(y,\omega)$ and $z_t(\omega):=
z_t(y,\omega).$ Suppose first that for $\omega \in \Omega,
Y_t(\omega) \rightarrow 0$ weakly in ${\cal S}^{\prime} $ and assume
that $z_t(\omega) \nrightarrow \infty. $ Since the neighborhoods of
$\infty$ in ${\mathbb R}^d$ are complements of compact sets , if
$z_t(\omega) \nrightarrow \infty, $ then there exists a ball
$B(0,r)$ of radius $r$ around zero and a sequence $t_n \uparrow
\eta(\omega)$ such that $z_{t_n}(\omega) \in B(0,r)$ for all $n \geq
1.$ The compactness of $B(0,r)$ implies the existence of a
subsequence of $(t_n)$, denoted again by $t_n,$ and $z \in B(0,r)$
such that $z_{t_n}(\omega) \rightarrow z.$ The continuity of the
translations and the weak convergence of $Y_{t_n}(\omega)$ to zero
now forces $\tau_z(y) = 0.$ This implies $y = 0,$ a contradiction.

For the converse suppose first that $y \in L^2({\mathbb R}^d).$ If
$z_t(\omega) \rightarrow \infty $ as $t \rightarrow \eta(\omega)$
then for $\phi \neq 0 \in {\cal S},$\Bea \langle \phi ,Y_t(\omega)
\rangle  &=& \int\limits_{\{x : |x| < n\}}y(x)\phi(x -
z_t(\omega))~dx + \int\limits_{\{x : |x| \geq n\}}y(x)\phi(x -
z_t(\omega))~dx \Eea Since  $y \in L^2({\mathbb R}^d),$ the second
integral can be made small, independent of $t$ by choosing $n$
large. For the first integral, we can choose $t$ sufficiently close
to $\eta(\omega)$ so that $\phi(x-z_t(\omega))$ is small, uniformly
for $|x| \leq n, $ proving the case when $y \in L^2({\mathbb R}^d).$

Suppose now that the second case holds i.e. $y$ has compact support.
Let support $y \subseteq K$ and let $N=$ order$(y)+2d$. Then there
exist continuous functions $g_\alpha , |\alpha |\leq N,
\mbox{support}~ g_\alpha \subseteq V$ where $V$ is an open set
having compact closure, containing $K$, such that \Bea y
=\sum\limits_{|\alpha | \leq N} \partial^\alpha g_\alpha . \Eea See
\cite{FT}, Theorem 24.5, Corollary 3. Then for $\phi \neq 0 \in
{\cal S},$\Bea \langle \phi ,Y_t(\omega) \rangle &=&
\sum\limits_{|\alpha | \leq N}\langle \phi ,
\tau_{z_t(\omega)}(\partial^\alpha g_\alpha) \rangle \\
&=& \sum\limits_{|\alpha | \leq N}(-1)^{|\alpha|}\langle
\partial^\alpha\phi , \tau_{z_t(\omega)}(g_\alpha)
\rangle.\Eea The same arguments as in the first case applied to each
of the terms in the above sum  will now also prove the second case.
$\hfill{\Box}$
\begin{Remark} \rm {It is easy to see that $z_t(y,\omega) \rightarrow \infty $
does not in general imply that $Y_t(y,\omega) \rightarrow 0$ weakly
in ${\cal S}^{\prime} $. For example if $d=1$ and $y = c$ , a non
zero constant , then $ y \in {\cal S}_{-p}$ for some $p > 0.$ If
$\sigma , b$ are integrable functions with non zero integrals over
${\mathbb R}$ then it is easy to see that $\eta = \infty$ and
$z_t(y,\omega) \rightarrow \infty$ almost surely on the one hand ,
while on the other $Y_t(y,\omega) = c$ for all $t \geq 0.$ }
\end{Remark}
\section{$(Y_t(y))$ as a Markov Process on $\hat{{\cal S}}_p.$}

In this section, we study the Markov properties of the solutions of
equation (3.4) viz. $(Y_t(y)).$ For this purpose, it is essential to
obtain a version of $(Y_t(y))$ which is jointly measurable in $y$ as
well. It is of course no accident and certainly not unreasonable
that the Markov properties of $Y_t(y)$ derive from that of the
diffusion $(X(x,y,t))$ generated by equation(2.3).In Proposition 4.1
below, we obtain a version of $((X(x,y,t)))$ which is also jointly
measurable in $x$ and $y.$ Let $\sigma_{ij},b_i \in {\cal S}_{-p},
i,j = 1,\cdots ,d$ and $y \in {\cal S}_p.$ Let
$\bar{\sigma}_{ij}(x,y) = \langle \sigma_{ij} , \tau_x y \rangle $
and $\bar{b}_i(x,y) = \langle b_i , \tau_x y \rangle $ be locally
Lipschitz functions on ${\mathbb R}^d$ (earlier denoted by
$\bar{\sigma}_{ij}(x),\bar{b}_i(x)$).
 Then by Theorem 2.1, equation(2.3)has a solution for each $x \in {\mathbb R}^d.$
 We will denote the Borel $\sigma$ field on $\hat{\mathbb R}^d $ by ${\cal B}_d.$

\begin{Proposition} Let  $\left(\Omega , {\mathcal F}, \left\{{\mathcal
F}_t\right\}_{t\geq 0}, P\right)$ be a filtered probability  space
satisfying the {\it usual conditions} and $(B_t)$ be a standard
${\cal F}_t$ Brownian motion on it. Then,there exists a map $X :
\hat{{\mathbb R}}^d \times \hat{{\cal S}}_p \times
[0,\infty)\times\Omega \rightarrow \hat{{\mathbb R}}^d $ which is
${\cal B}_d\otimes {\cal B}(\hat{{\cal S}}_p)\otimes {\cal
B}[0,\infty)\otimes {\cal F}$ measurable and such that for each
$(x,y) \in {\mathbb R}^d \times {\cal S}_p, (X(x,y,t))$ is a
solution of equation (2.3) given by Theorem 2.1.

\end{Proposition}

{\bf Proof:}  We note that $\bar{\sigma}(x,y)$ and $\bar{b}(x,y)$
are jointly measurable in $(x,y).$ Choose for each $n \geq 1,$ and
$i,j = 1,\cdots d,$ Lipshitz functions in $x$ $\sigma_{ij}^n(x,y),
b_i^n(x,y) $ , measurable in $(x,y)$such that $\sigma_{ij}^n(x,y) =
\bar{\sigma}_{ij}(x,y),~ b_i^n(x,y) = \bar{b}_i(x,y), |x| \leq n.$
Consider equation(2.3)  with $\bar{\sigma}_{ij},~\bar{b}_i$ replaced
by $\sigma_{ij}^n(x),{b}_i^n(x)$ viz. \bea dX_t =
\sigma^n(X_t,y).dB_t + b^n(X_t,y) dt ~~~; ~~~X_0 = x. \eea Let for
$n \geq 1, k \geq 1,$$(X^{n,k}(x,y,t))$ be defined iteratively (in
vector form) by $$ X^{n,k}(x,y,t) := x + \int_0^t
\sigma^n(X^{n,k-1}(x,y,s),y).dB_s + \int_0^t
b^n(X^{n,k-1}(x,y,s),y)~ds $$  with $X^{n,0}(x,y,t) := x$ for all $t
\geq 0.$ It is a well known property of stochastic integrals that
they are measurable with respect to a parameter when the integrands
are measurable with respect to the same parameter (see \cite{OK},
Theorem 17.25). Using an inductive argument and the fact that
$\sigma_{ij}^n(x,y),b_i^n(x,y)$ are jointly measurable in $(x,y)$ it
follows that the integrals in the right hand side of the above
equation have jointly measurable versions in $(x,y,t,\omega)$ and
the same follows for $X^{n,k}(x,y,t,\omega).$ Then by the method of
successive approximations , we get solutions
$(X^n(x,y,t,\omega))$,of equation(4.6), where $X^n(x,y,t,\omega) =
\lim\limits_{k \rightarrow \infty} X^{n,k}(x,y,t,\omega)$ , jointly
measurable in $(x,y,t,\omega)$ and for each $(x,y)$, progressively
measurable in $(t,\omega)$(see \cite{IW},Chapter IV, Theorem 3.1).
Note that for each $(x,y), (X^n(x,y,t,\omega))$ also solves
equation(2.3) upto the random time $\eta^n(x,y,\omega)$ defined as
\Bea \eta^n(x,y,\omega)&:=& \rm inf \{ t \geq 0 : \|
X^n(x,y,t,\omega) \| \geq n \}  \\ &= & \rm inf \{ t \geq 0 : \|
X(x,y,t,\omega) \| \geq n \}.\Eea It is easy to see that
$\eta^n(x,y,\omega)$ is jointly measurable in $(x,y,\omega).$
Further,if we denote by $\eta(x,y,\omega)$ the explosion time for
the solution $(X(x,y,t))$ of equation (2.3) starting at $x \in
{\mathbb R}^d,$ then $\eta(x,y,\omega) := \lim\limits_{n \rightarrow
\infty }\eta^n(x,y,\omega).$ As a consequence, $\eta(x,y,\omega)$ is
a measurable function of $(x,y,\omega) \in {\mathbb R}^d\times {\cal
S}_p \times \Omega.$ The map $X$ may now be  defined as
$$X(x,y,t,\omega) := \lim\limits_{n \rightarrow \infty}
X^n(x,y,t,\omega)I_{\{t < \eta(x,y,\omega)\}}(x,y,t,\omega) + \infty
I_{\{t \geq \eta(x,y,\omega)\}}(x,y,t,\omega).$$ Clearly for each $
(x,y) \in {\mathbb R}^d \times {\cal S}_p$ we have , by the
uniqueness of solutions, $X(x,y,t\wedge \eta^n) = X^n(x,y,t \wedge
\eta^n) .$ We define $X(x,y,t,\omega) = \infty $ for all
$(t,\omega)$ if $x = \infty$ or $ y = \delta.$ $\hfill{\Box}.$

We define the transition probability function $\bar{P}(x,y,t,A)$ of
the diffusion $(X(x,y,t))$ in the usual way : For $0 \leq t \leq
\infty, x \in {\mathbb R}^d, y \in {\cal S}_p$ and $ A \in {\cal
B}_d,$ \Bea \bar{P}(x,y,t,A) &:=& P\{\omega: X(x,y,t,\omega) \in A
\} \\&=& P\{ X(x,y,t) \in A\cap{\mathbb R}^d, t < \eta(x,y)\} +
I_{A}(\infty)P\{ t \geq \eta(x,y)\}. \Eea Note that $t = \infty$ is
included in the definition of $\bar{P}(x,y,t,A)$ by taking
$X_{\infty} = \infty.$We take $\bar{P}(x,y,t,A) := I_A(\infty), t
\geq 0,$ if $ x = \infty .$ We can define an induced transition
probability on ${\cal S}_p$ as follows : First we extend the map
$\tau_x(y).$ We define $\tau_{\infty}(y) := \delta , y \in
\hat{{\cal S}}_p$ and $\tau_x(\delta) := \delta.$ Thus $\tau_x :
\hat{{\cal S}}_p \rightarrow \hat{{\cal S}}_p , x \in \hat{{\mathbb
R}^d}.$ For $y \in {\cal S}_p, 0 \leq t \leq \infty ,$ and $B \in
{\cal B}(\hat{{\cal S}}_p)$ define
$$P(y,t, B) :=  \bar{P}(0,y,t,
\tau^{-1}_{.}(y)(B))$$ where $ \tau^{-1}_{.}(y)(B) = \{z:
\tau_{z}(y) \in B \}. $ We take $P(y,t,B) := I_B(\delta), t \geq 0,$
if $y = \delta$ and define $Y_t(y) = \delta$ if $t = \infty$ or $y =
\delta.$ We note that because $X(0,y,t) = z_t(y), t <
\eta(0,y)$ we have  \Bea P(y,t,B) &=& P\{X(0,y,t) \in \tau^{-1}(y)(B)\} \\
&=& P\{ X(0,y,t) \in \tau^{-1}(y)(B), t < \eta(0,y)\} \\&& + P\{t
\geq \eta(0,y)\}I_{\tau^{-1}(y)(B)}(\infty) \\  &=& P\{
\tau_{z_t(y)}(y) \in B, t < \eta(0,y)\} + P\{t \geq
\eta(0,y)\}I_{B}(\delta)\\ &=& P\{Y_t(y) \in B \} \Eea The strong
Markov property for $(Y_t(y))$ is now a simple consequence of that
for the process $(X(x,y,t)).$

\begin{Proposition} Let $y \in {\cal S}_p$ and let $T$ be an ${\cal F}_t $
 stopping time . Then for $ 0 \leq s \leq \infty ,$ and $B \in {\cal B}(\hat{{\cal S}}_p),$ we have ,
almost surely , $$ P\{ Y_{s +T}(y)\in B | {\cal F}_T \} = P(s,
Y_T(y), B)$$
\end{Proposition}
{\bf Proof :} Since the result holds trivially at $s = \infty,$ we
assume $s < \infty.$ We have, using the strong Markov property of
the process $(X(x,y,t))$ \Bea P\{ Y_{s +T}(y)\in B
| {\cal F}_T \} &=& P\{ \tau_{z_{s +T}(y)}(y)\in B | {\cal F}_T \} \\
&=& P\{ z_{s +T}(y)\in \tau_{.}^{-1}(y)(B) | {\cal F}_T \}  \\ &=&
P\{ X(0,y, s +T)\in \tau_{.}^{-1}(y)(B) | {\cal F}_T \} \\ &= &
\bar{P}(X(0,y,T),y,s,\tau_{.}^{-1}(y)(B) )\Eea On the other hand,
\Bea \bar{P}(x,y,s, \tau^{-1}(y)(B)) &=& P\{X(x,y,s) \in
\tau^{-1}(y)(B)\} \\ &=& P\{x + z_s(\tau_x(y)) \in \tau^{-1}(y)(B),
s < \eta(x,y)\} \\ && + P\{ s \geq
\eta(x,y)\}I_{\tau^{-1}(y)(B)}(\infty)\\ &=& P\{\tau_{x +
z_s(\tau_x(y))}(y) \in B, s < \eta(x,y)\} + P\{ s \geq
\eta(x,y)\}I_{B}(\delta)\\ &=& P\{\tau_{ z_s(\tau_x(y))}(\tau_x(y))
\in B, s < \eta(x,y)\} + P\{ s \geq \eta(x,y)\}I_{B}(\delta)\\
&=& P(\tau_x(y),s,B)\Eea where in the last equality we have made use
of the observation made preceeding the statement of the proposition
and the fact that $Y_s(\tau_x(y)) =
\tau_{z_{s}(\tau_x(y))}(\tau_{x}(y))$. Hence \Bea P\{ Y_{s +T}(y)\in
B | {\cal F}_T \} &=& \bar{P}(X(0,y,T),y,s,\tau_{.}^{-1}(y)(B) ) \\
&=& P(\tau_x(y),s,B)|_{x= X(0,y,T)} \\ &=& P(\tau_x(y),s,B)|_{x=
z_{T}(y)} \\ &=& P(Y_T(y),s,B)\Eea This completes the proof of the
Proposition.$\hfill{\Box}$

Let ${\mathbb B}_p = \{ f : \hat{{\cal S}}_p \rightarrow {\mathbb R}
; f \rm {~bounded~and~ measurable}, f(\delta) = 0 \}.$ Let
$\|.\|_{p,\infty}$ denote the norm on ${\mathbb B}_p$ given by
$\|f\|_{p,\infty} := \sup\limits_{y \in {\cal S}_p}|f(y)|.$ Then
$({\mathbb
 B}_p ,\|.\|_{p,\infty})$ is a Banach space.  Let $(T_t)_{0 \leq t < \infty}$
denote the linear operators $T_t : {\mathbb
 B}_p \rightarrow {\mathbb B}_p $ given by $T_tf(y) := E(f(Y_t(y))),$ for $ y \in {\cal S}_p$
and $f \in {\mathbb B}_p.$
\begin{Corollary} $(T_t)_{t \geq 0}$ is a semi-group of linear
operators on ${\mathbb B}_p$ with $T_0 = Identity ,~ T_t1 = 1$ and
$T_tf \geq 0$ whenever $f \geq 0.$
\end{Corollary}
Let $({\cal L}, Dom({\cal L}))$ denote the infinitesimal generator
of $(T_t).$ Recall that $f \in {\mathbb B}_p$ belongs to $Dom({\cal
L})$ iff the limit of $\frac{1}{t}(T_tf - f)$ as $t$ tends to zero
exists in ${\cal B}_p$ and further ${\cal L}(f) := \lim\limits_{t
\rightarrow 0}\frac{1}{t}(T_tf - f).$ We shall denote by
$(\bar{T}_t^y)$ and $\bar{{\cal L}}^y$ the semi-group and
infinitesimal generator, respectively on ${\mathbb B}_d$ and $
Dom(\bar{{\cal L}}^y) \subset {\mathbb B}_d,$ associated with the
diffusion $(X(x,y,t))$ generated by $(\bar{\sigma}(x,y))$ and
$\bar{b}(x,y).$ Here ${\mathbb B}_d$ is the Banach space of bounded
measurable functions on $\hat{{\mathbb R}}^d$ endowed with the
supremum norm.

For $y \in {\cal S}_p,$ let $C(y) \subset {\cal S}_p$ be defined by
$C(y) := \{y^{\prime}  \in {\cal S}_p : y^{\prime} = \tau_x(y), x
\in {\mathbb R}^d \}.$ $\hat{C}(y) := C(y)\cup \{\delta\}$.Then
observe that ${\cal S}_p = \bigcup\limits_{y \in {\cal S}_p}C(y)$
and for $y_1 \neq y_2$ either $C(y_1)\bigcap C(y_2) = \phi$ or
$C(y_1) = C(y_2).$ We shall consider $C(y)$ as a measurable space
with the $\sigma-$field induced by ${\cal B}(\hat{{\cal S}}_p).$ Let
$${\mathbb B}_p(y) := \{f:\hat{C}(y) \rightarrow {\mathbb R}; f
\rm{~bounded~and~measurable},f(\delta) = 0 \}.$$ Since $P\{ Y_t(y)
\in C(y^{\prime}) \rm{~for~some~}t \geq 0 \} =0 $ if $C(y)\bigcap
C(y^{\prime}) = \phi $ we have (the restriction) $T_t|_{{\mathbb
B}_p(y)} : {\mathbb B}_p(y) \rightarrow {\mathbb B}_p(y)$ is a
semi-group on $({\mathbb B}_p(y),\|.\|_{p,\infty}),$ for every $y
\in {\cal S}_p.$ We shall continue to denote the restrictions of
$T_t$ and ${\cal L}$ to ${\mathbb B}_p(y)$ and ${\mathbb
B}_p(y)\bigcap Dom({\cal L})$ respectively by $(T_t)$ and ${\cal L}$
or by $T_t^y$ and ${\cal L}^y$ if there is a risk of confusion.

  Let $F: {\mathbb R}^{d+n} \rightarrow {\mathbb R},
 F \in C^{\infty}({\mathbb R}^{d+n})$ such that $\rm {support}(F) \subset K \times {\mathbb R}^n$ for some compact $K \subset {\mathbb R}^d. $
 Let $\varphi_i \in {\cal S}, i = 1,\cdots d.$ Fix $y \in {\cal S}_p.$Define $f \equiv f^y : {\mathbb R}^d \times {\cal S}_p \rightarrow {\mathbb R}$ as follows :
If $y^{\prime} = \tau_x(y)$ define $$ f^{y}(x,y^{\prime}) :=
F(x,\langle
\varphi_1,\tau_x(y)\rangle,\cdots,\langle\varphi_n,\tau_x(y)\rangle).$$
If $y^{\prime} \notin C(y)$ we define $f^y(x,y^{\prime}) = 0.$
Define $\bar{f}^y(x) := f^y(x,\tau_x(y)) : {\mathbb R}^d \rightarrow
{\mathbb R}$ and $\bar{f}^y(\infty) = 0.$ Let $\bar{\sigma}^t(x,y)$
 denote the transpose of the matrix given by $\bar{\sigma}(x,y) :=
(\bar{\sigma}_{ij}(x,y)).$ We will denote by $\bar{L}^y$ the second
order differential operator in the variable $x$ given as
$$ \bar{L}^y := \sum\limits_{i,j =1}^d
(\bar{\sigma}(x,y)\bar{\sigma}^t(x,y))_{ij}~\partial_{ij}^2 +
\sum\limits_{i=1}^d \bar{b}_i(x,y)\partial_i .$$

\begin{Proposition} Let $y \in {\cal S}_p.$ Then $\bar{f}^y(x) \in C^{\infty}_K({\mathbb R}^d) \subset Dom(\bar{L})$.
 Consequently,for $x \in {\mathbb R}^d,$
 $f^y(x,.) \in  Dom({\cal L})$ and $$ {\cal L}f^y(x,\tau_x(y)) = {\cal L}^yf^y(x,\tau_x(y)) = \bar{{\cal L}}^y\bar{f}^y(x) = \bar{L}^y\bar{f}^y(x) .$$
  ${\cal L}f^y(x,y^{\prime}) = 0 $ if $y^{\prime} \notin C(y).$
\end{Proposition}
{\bf Proof:} It is clear that $\bar{f}^y \in C^{\infty}_K.$ Further,
using the compactness of $K$ and Ito's formula it can be shown that
$C^{\infty}_K \subset Dom(\bar{L}).$ It is easily seen from the
definitions that $T_tf^y(x,\tau_x(y)) = \bar{T}_t^y\bar{f}^y(x)$ and
$T_tf^y(x,y^{\prime}) = 0, y^{\prime} \notin C(y).$ In particular,
$$\frac{1}{t}\{T_tf^y(x,y^{\prime}) - f^y(x,y^{\prime})\} =
I_{\{\tau_x(y): x \in {\mathbb
R}^d\}}(y^{\prime})~\frac{1}{t}\{\bar{T}_t^y\bar{f}^y(x) -
\bar{f}^y(x)\} .$$ The result follows.$\hfill{\Box}$

\begin{Remark} \rm {Let $ F(x,z) := g(x)z ,$ where $ g \in C^{\infty}_K , K =
B(0,r),$ the ball of radius $r$ centred at $0,$ and $g(x) = 1 , x
\in B(0,r_1) $ for some $r_1 < r.$ Then for $y \in {\cal S}_p$
$f^y(x,y^{\prime}) = g(x)\langle \phi,y^{\prime}\rangle, \phi \in
{\cal S}, y^{\prime} \in C(y) ; f^y(x,y^{\prime}) = 0 , y^{\prime}
\notin C(y) $. It is easy to see that
$${\cal L}f^y(0,y) = \bar{L}^y\bar{f}^y(0) = \langle \phi, L(y) \rangle  $$
where $L$ is the operator in equation(3.4).}
\end{Remark}

\begin{Remark} \rm {Assume $y \in {\cal S}_p$ is such that  $\tau_xy \neq y$
for any $x \neq 0.$ Let $j: \hat{C}(y)\rightarrow \hat{{\mathbb
R}}^d $ be defined by $ j(\tau_x(y)) := x.$ Then $j$ is one-one and
onto and we provide $\hat{C}(y)$ with a topology and corresponding
Borel structure ${\cal B}_y $ that makes $j: (\hat{C}(y),{\cal B}_y)
\rightarrow (\hat{{\mathbb R}}^d, {\cal B}_d)$ a Borel isomorphism
with inverse $\tau_.(y):\hat{{\mathbb R}}^d \rightarrow \hat{C}(y)$.
We can extend $j$ as a map $ j: {\mathbb B}_p(y) \rightarrow
{\mathbb B}_d $ such that $\bar{f} \in {\mathbb B}_d $ iff $ f \in
{\mathbb B_p(y)},$ where $f = \bar{f}\circ j.$ This extends to semi
groups viz. $ T_t^y = \bar{T}_t^y \circ j.$ In other words, the
Markov process $(Y_t(y)) \equiv ( \tau_{z_t(y)}(y))$ on the state
space $\hat{C}(y)$ is `isomorphic' to the Markov process
$(X(x,y,t)))$ on $\hat{{\mathbb R}}^d.$ Note however that the
topology on $\hat{C}(y)$ is that of $\hat{{\mathbb R}}^d $ and is
different from the one induced from $\hat{{\cal S}}_p.$ Proposition
3.12 is a reflection of the same phenomenon.}
\end{Remark}

\section{ A Non Linear Evolution Equation .}
 In this section  we derive a non
 linear evolution equation associated with
the operator $L$ with the initial condition $y \in S_p$ viz. \bea
\partial_t \psi (t,y) = \psi (t,L(y))~~~;~~~
\psi (0,y) = y. \eea  We construct solutions of (5.7) via what maybe
called non linear convolutions, that we define below.  This is also
closely related to the notion of stochastic representation of
solutions to evolution equations of the type (4.6). See
\cite{RT1},\cite{RT2}. The solutions of equation(5.7)are also
related to the solutions of the forward equation for the diffusion
$(X(x,y,t))$(see below, equation(5.8), Remark 5.5.) We recall from
\cite{RT1} that $\tau_x :{\cal S}_p \rightarrow {\cal S}_p$ are
bounded linear operators for all $p \in {\mathbb R}$.

\begin{Definition}\rm{Let $p \in {\mathbb R}$ and let $ q \leq p$. Suppose
$h : {\cal S}_p \rightarrow {\cal S}_q$ and $f : {\mathbb R}^d
\rightarrow {\mathbb R}$ be  Borel measurable maps. For $ y \in
{\cal S}_p$, the convolution $h(y) \circ f$ is defined to be the
element of ${\cal S}_q$ given by the Bochner integral $ h(y) \circ f
:= \int\limits_{{\mathbb R}^d} h(\tau_xy) f(x)~dx$ provided the
integral exists i.e. provided $\int\limits_{{\mathbb R}^d}
\|h(\tau_x(y))\|_q|f(x)|~dx < \infty$.More generally, let $\mu$ be a
finite measure on the Borel sigma field of ${\mathbb R}^d$ and
$h(y)$ be as above. The convolution $h(y) \circ \mu$ is defined as $
h(y) \circ \mu := \int\limits_{{\mathbb R}^d} h(\tau_xy) \mu(dx) $
provided $\int\limits_{{\mathbb R}^d}\|h(\tau_x(y))\|_q~\mu(dx) <
\infty$. }
\end{Definition}
\begin{Remark}\rm { Our notation is a compromise between two contrasting
interpretations of the above definition. We could interpret the
above definition as an extension of the notion of convolution of a
functions $h : {\mathbb R}^d \rightarrow {\mathbb R}$ and a finite
measure $\mu$ on ${\mathbb R}^d,$ to that of convolution of $\mu$
and a map $h : {\cal S}_p \rightarrow {\cal S}_q.$ The notation then
would be $h \circ \mu (y),$ where $h\circ \mu :{\cal S}_p
\rightarrow {\cal S}_q.$ The definition also affords an
interpretation as an extension of the notion of convolution of a
tempered distribution $y$ and a measure $\mu$ via the map $h.$ We
may view it then as a non linear convolution between $y$ and $\mu$
where the nonlinearity arises because of the map $h.$ An appropriate
notation then could be $y \circ_h f $ or $y \star_h f.$} In any
case, our definition of convolution reduces to the usual convolution
between two distributions $y$ and $\mu$, if we take $q = p$ and
$h(y) = y$ in Definition $5.1.$ We will use the notation $y \star
\mu$ for the usual convolution between the distribution $y$ and the
measure $\mu.$ We further note that when $h$ is non linear as in
definition (5.1), $h(y)\circ \mu$ is , in general, different from
the ordinary convolution $h(y)\star \mu$, between the distribution
$h(y)$ and the measure $\mu.$. In our application of the notion of
convolution to construct solutions of equation(5.7) however, there
is an additional source of non-linearity viz. $\mu$ would also
depend on $y.$
\end{Remark}

\begin{Definition} Let $y \in {\cal S}_p.$ We say that a continuous map
$\psi(.,y) : [0,\infty) \rightarrow {\cal S}_p ,$ is a solution by
convolution of the initial value problem (5.7) iff there exists
kernels $\mu(t,dx) $ on $[0,\infty)\times{\mathbb R}^d$ such that
\begin{enumerate} \item  $\int\limits_{{\mathbb R}^d} \|\tau_x(y)\|_p~ \mu(t,dx) < \infty,
t \geq 0 ,$\rm {~and~} $\int\limits_{{\mathbb R}^d} \|\phi(x)
\tau_x(y)\|_p \mu(t,dx)~ < \infty ,\\ \rm{for~  all~ } t \geq
0,\rm{and ~ for ~all}~ \phi,\rm{where}~ \phi(x) = (\langle \sigma
,\tau_x(y)\rangle \langle\sigma,\tau_x(y)\rangle^t)_{ij} , i,j=
1,\cdots d,\rm{or}~ \phi(x)= \langle b_{i},\tau_x(y)\rangle, i = 1,
\cdots,d. $ In particular if $L : {\cal S}_p \rightarrow {\cal
S}_{p-1}$ is as in equation (3.4), then the convolution  $L(y)\circ
\mu(t,.) $ exists in ${\cal S}_{p-1}$ for all $t \geq 0.$
\item $\psi(t,y) = y \circ \mu_t , t \geq 0$
\item $\psi(t,y)$ is continuously differentiable for $t \in (0,\infty)$ and we have
$$ \partial_t\psi(t,y)= \partial_t(y \circ \mu_t) = L(y)\circ \mu(t,.)~~~;~~~ \psi(0,y) = y \circ \mu_0= y.$$

\end{enumerate}
\end{Definition}

Recall from Section 4 ,the transition probability measure for the
solutions $(Y_t(y))$ of equation(3.4) viz. $P(y,t,B) $ and the
transition probability measure for the process $(z_t(y))$ given by
Corollary 3.8 viz. $\bar{P}(0,y,t,A).$ We recall that $P(y,t,B) =
\bar{P}(0,y,t,\tau_.^{-1}(y)(B)).$ The following theorem constructs
the solutions of the initial value problem  equation(5.7) via a
stochastic representation.

\begin{Theorem} Let $\sigma_{ij},b_i, y$ be as in Theorem 3.4 and
let $\{(Y_t(y)), \eta \}$ be the unique $S_p$ valued solution of
equation (3.4) given by Theorem 3.11. Let $\bar{\sigma}_{ij}(x,y) :=
\langle \sigma_{ij}, \tau_x(y)\rangle , \bar{b}_i(x,y) := \langle
b_i, \tau_x(y)\rangle, i,j = 1,\cdots,d $ and suppose that for fixed
$y,$ these are bounded and continuous functions of $x.$ Let
$(z_t(y))$ be the unique solution of equation (2.3), as in Theorem
3.11. Then $\psi(t,y) := E(Y_t(y)), t \geq 0$ defines an ${\cal
S}_p$ valued continuous map that solves the initial value problem
(5.7) by convolution .
\end{Theorem}

{\bf Proof : } We note that under the assumptions on
$\bar{\sigma}_{ij}(.,y),\bar{b}_i(.,y),$ $(z_t(y))$ has moments of
all orders and further for all $t \geq 0,$
$$ \sup\limits_{s \leq t}E|z_s(y)|^{k} < C(t,y,k) < \infty ~~~ k = 1,2,\cdots .$$
In particular, $\eta = \infty , a.s.$
 Note that $\psi(t,y)$ is well defined : for $z \in {\mathbb R}^d$,
 $$\|\tau_z y\|_p ~~\leq~~ ~\|y\|_p P(|z|) $$
 where $P(x)$ is a polynomial in $x \in {\mathbb R}$ of degree $2|p|+1$
 (see \cite{RT1})  .
 It follows from our assumption on the moments of $(z_t)$
that $$E\|Y_t\|_p ~~\leq~~ ~\|y\|_p EP(|z_t|) ~~<~~ \infty$$ for
every $t \geq 0$. In particular $\psi(t,y) : =  EY_t(y) $ exists as
a Bochner integral in $S_p$. Further it is clear that  $\psi(t,y) =
y \circ \bar{P}(0,y,t,.)$ The theorem is proved by taking
expectations in equation (3.5) satisfied by $(Y_t(y)).$ We first
show that this is indeed a legitimate operation.

Using the fact that the moments of $z_s := z_s(y)$ are finite, we
have for each $t \geq 0$, \Bea E\|\int_0^tL(\tau_{z_s}y)~ds\|_{p-1}
&\leq&
 E\int_0^t\|L(\tau_{z_s}y)\|_{p-1}~ds  \\ &\leq&  C^{\prime} E\int_0^t
 \{\|\tau_{z_s}y\|_p^3 + \|\tau_{z_s}y\|_p^2 \}~ds \\ &\leq &
 C E\int_0^t\{ P(|z_s|)^3 + P(|z_s|)^2 \}~ds \\ & < &  \infty
\Eea where $C^{\prime} =
C^{\prime}(d,\|\sigma_{i,j}\|_{-p},\|b_i\|_{-p}, i,j = 1 \cdots d)$
and $C = C(d,\|\sigma_{i,j}\|_{-p},\|b_i\|_{-p},\\ i,j = 1 \cdots
d,\|y\|_p)$ are positive constants depending on the indicated
quantities. A similar calculation verifies that for each $t \geq 0$
\Bea
 E \| \sum\limits_{j=1}^d\int\limits_0^t A_j(\tau_{z_s}y)) ~dB_s^j \|^2_{p-1}& \leq
& \sum\limits_{j=1}^d E\int\limits_0^t\| A_j(\tau_{z_s}y)\|^2_{p-1}
~ds
\\ &\leq & C E\int_0^t\{ P(|z_s|)^4 \}~ds \\ &<& \infty \Eea We can
thus take expectations in equation (3.5) to get
 \Bea
\psi(t,y) = EY_t &=&  y +E \int\limits_0^t L(\tau_{z_s}y)~ds \\
&=& y +\int\limits_0^t ~\int\limits_{\mathbb R^d} L(\tau_x y)~
\bar{P}(0,y,s,dx)~ds \\
&=& y + \int\limits_0^t ~L(y) \circ \bar{P}(0,y,s,.)~ds \\
&=& y + \int_0^t E(L(Y_s(y)))~ds \Eea That $\psi(t,y)$ is
continuously differentiable and satisfies equation (5.7)now follows
from the above equation and the continuity of $E(L(Y_s(y))).$ This
completes the proof of Theorem 5.4. $\hfill{\Box}$

For $\phi, y \in {\cal S}^{\prime}$ such that the products
$\bar{\sigma}_{ij}(.,y) \phi,~\bar{b}_i(.,y)\phi, i,j = 1,\cdots d$
are tempered distributions in the variable $x$, we define the
operator $\bar{L}^{*,y} $ as follows :
$$\bar{L}^{*,y} \phi  := \frac{1}{2}~\sum\limits^d_{i,j=1}
\partial^2_{ij}((\bar{\sigma}(.,y)\bar{\sigma}^t(.,y))_{ij} \phi) -\sum\limits^d_{i=1}~ \partial_i (\bar{b}_i(.,y)\phi).$$

We note that $\bar{L}^{*,y}$ is the formal adjoint of the second
order differential operator $\bar{L}$ associated with the diffusion
with coefficients $\bar{\sigma}_{ij},\bar{b}_i , i,j = 1, \cdots
,d.$ The initial value problem (5.7)is closely connected with
solutions of the forward Kolmogorov equation for $\bar{L}$ viz. \bea
\partial_t \psi_t = \bar{L}^{*,y}\psi_t~~~~;~~~~ \psi_0 =\psi  \eea
When $\psi$ is a distribution with compact support and
$\bar{\sigma}_{ij}, \bar{b}_i$ are smooth , solutions to the above
equation maybe obtained by convolution with the transition
probability measure $\bar{P}(0,y,t,.)$ (see \cite{RT2}, Theorem
4.5.). We extend that result in Theorem 5.6 when the coefficients
are only bounded and continuous.

\begin{Remark} \rm{ To see the connection between solutions of equation (5.7)
and solutions of equation(5.8), consider the following integrated
version of equation(5.8) with $\psi_t = \bar{P}(0,y,t,.)$ viz. \bea
\bar{P}(0,y,t,.) = \delta_0 + \int_0^t \bar{L}^*\bar{P}(0,y,s,.)
~ds\eea

 By convolving with $y \in {\cal S}_p $ and using the relation
$$y\star\bar{L}^{*,y}\bar{P}(0,y,t) = L(y)\circ\bar{P}(0,y,t,.)$$
we get , $$ y \circ \bar{P}(0,y,t,.) = y + \int_0^t
L(y)\circ\bar{P}(0,y,s,.)~ds $$ which is equivalent to equation
(5.7). On the other hand we can Fourier transform the above to get
back (5.9): Suppose $y$ has compact support. For a tempered
distribution $\phi,$ the Fourier transform of $\phi$ is denoted by
$\hat{\phi}.$ For each $\xi \in {\mathbb R}^d,$\Bea
\hat{y}(\xi)\bar{P}(0,y,t,.\hat{)}(\xi) &=& \hat{y}(\xi) + \int_0^t
(L(y)\circ \bar{P}(0,y,s,.)\hat{)}(\xi)~ds
\\ &=& \hat{y}(\xi) + \int_0^t
(\bar{L}^{*}\bar{P}(0,y,s,.)\hat{)}(\xi)\hat{y}(\xi)ds\Eea Since $y$
has compact support, $\hat{y}(\xi),$ is by the Paley-Wiener theorem
an entire function and hence we may cancel off $\hat{y}(\xi)$ in the
above equation to get, almost surely with respect to Lebesgue
measure on ${\mathbb R}^d,$ \Bea \bar{P}(0,y,t,.)\hat{)}(\xi) = 1 +
\int_0^t (\bar{L}^{*} \bar{P}(0,y,s,.)\hat{)}(\xi)ds \Eea }
Inverting the Fourier transform , we get back equation (5.9).
\end{Remark}

\begin{Theorem} Suppose that $y \in {\cal S}_p$, $\sigma_{ij}, b_i, \in {\cal S}_{-p}, i,j = 1,\cdots d$
and in addition are such that $\bar{\sigma}_{ij}(.,y),
\bar{b}_i(.,y)$ are bounded continuous functions. Let $q >
\frac{d}{4}.$ Then for any $x \in {\mathbb R}^d,$the map $t
\rightarrow \bar{P}(x,y,t,.) : [0,\infty) \rightarrow {\cal S}_{-q}$
is differentiable and satisfies the forward equation (5.8) with
$\psi = \delta_x.$

\end{Theorem}
{\bf Proof :} Let $(X(x,y,t))$ be the unique solution of (2.3). Then
by our assumptions on $\bar{\sigma}_{ij},\bar{b}_i$ , the moments of
all orders of $(X(x,y,t))$ exist and are locally bounded functions
of $t.$. Further since $\delta_z \in {\cal S}_{-q}$ if and only if
$q > \frac{d}{4}$(see \cite{RT2}), we have $\delta_{X(x,y,t)} \in
{\cal S}_{-q}.$ In particular, as in the proof of the previous
theorem, $E\|\delta_{X(x,y,t)}\|_{-q} < \infty.$ It follows that
$\bar{P}(x,y,t,.) = E\delta_{X(x,y,t)},$ where the right hand side
is an element of ${\cal S}_{-q}$ and the equality holds there. Using
the Ito formula in \cite{BR} we get, \bea \delta_{X(x,y,t)} &=&
\delta_x + \int_0^t \bar{L}(s)(\delta_{X(x,y,s)})~ds
 + \int_0^t \bar{A}(s)(\delta_{X(x,y,s)})\cdot dB_s~~
\eea where the  operator valued processes $(\bar{L}(s,\omega))$ and
$(\bar{A}(s,\omega))$ are defined for fixed $x $ and $y$ as follows.
$\bar{A}(s,\omega) := (\bar{A}_1(s,\omega),\cdots,
\bar{A}_d(s,\omega))$ and for $\phi \in {\cal S}^{\prime}$,
 $$\bar{A}_i(s,\omega)
\phi := -\sum\limits^d_{j=1}~  \bar{\sigma}_{ji}(X(x,y,s,\omega),y)
~\partial_j \phi.$$ Similarly \Bea\bar{L}(s,\omega)\phi &:=&
\frac{1}{2}~\sum\limits^d_{i,j=1} ( \bar{\sigma}(X(x,y,s,\omega),y)
 \bar{\sigma}^t(X(x,y,s,\omega),y))_{ij}~ \partial^2_{ij} \phi \\ &&
-\sum\limits^d_{i=1}~\bar{b}_i(X(x,y,s,\omega),y)~ \partial_i
\phi.\Eea As in the proof of Theorem 5.4, we can take expectations
in equation (5.10) to get \bea \bar{P}(x,y,t,.) = \delta_x +
\int_0^t E(\bar{L}(s)\delta_{X(x,y,s)})~ds. \eea Note that we have
the representation $ \bar{P}(x,y,t,.) = \int\limits_{{\mathbb R}^d}
\delta_{z}(.)\bar{P}(x,y,t,dz),$ where the right hand side is a
Bochner integral in ${\cal S}_{-q}.$ Given  a bounded continuous
function $\phi$ we can use this representation to define the product
$\phi\bar{P}(x,y,t,.)$ as an element of ${\cal S}_{-q}$ as follows :
$$ \phi\bar{P}(x,y,t,.) = \int\limits_{{\mathbb R}^d}
\phi(z)\delta_{z}(.)\bar{P}(x,y,t,dz)$$ where the right hand side is
a Bochner integral in ${\cal S}_{-q}.$ Hence
$\bar{L}^*\bar{P}(x,y,t,.)$ is a well defined tempered distribution
in ${\cal S}_{-q-1}.$ It is now easy to see , by acting on test
functions, that $\bar{L}^{*,y}\bar{P}(x,y,t,.) =
E(\bar{L}(t)\delta_{X(x,y,t)})$. In particular from (5.11) we get
\bea \bar{P}(x,y,t,.) = \delta_x + \int_0^t
\bar{L}^*\bar{P}(x,y,s,.) ~ds.\eea Since the moments of $X(x,y,t)$
are locally bounded functions of $t$, it follows using the dominated
convergence theorem that the integrand  in the right hand side in
equation (5.12) is continuous in $t.$ The desired conclusion now
follows from equation (5.12). $\hfill{\Box}$

{\bf Acknowledgements :} The author would like to thank an anonymous
referee for his careful reading of the manuscript and for his useful
comments and suggestions.

\end{document}